\documentclass[pdflatex,iicol,sn-mathphys-num]{sn-jnl} 

\usepackage{graphicx}%
\usepackage{multirow}%
\usepackage{amsmath,amssymb,amsfonts}%
\usepackage{amsthm}%
\usepackage{mathrsfs}%
\usepackage[title]{appendix}%
\usepackage{xcolor}%
\usepackage{textcomp}%
\usepackage{manyfoot}%
\usepackage{booktabs}%
\usepackage{algorithm}%
\usepackage{algorithmicx}%
\usepackage{algpseudocode}%
\usepackage{listings}%

\usepackage{subcaption}
\usepackage{csml}

\usepackage{verbatim}

\theoremstyle{thmstyleone}%
\theoremstyle{thmstyletwo}%

\theoremstyle{thmstylethree}%

\begin{document}

\title[Neural parametric representations for thin-shell shape optimisation]{Neural parametric representations for thin-shell shape optimisation}

\author*[1]{\fnm{Xiao} \sur{Xiao}}\email{xiao.xiao@sjtu.edu.cn}

\author[2]{\fnm{Fehmi} \sur{Cirak}}\email{f.cirak@eng.cam.ac.uk}

\affil*[1]{\orgdiv{School of Ocean and Civil Engineering}, \orgname{Shanghai Jiao Tong University}, \orgaddress{\city{Shanghai}, \postcode{200240}, \country{China}}}

\affil[2]{\orgdiv{Department of Engineering}, \orgname{University of Cambridge}, \orgaddress{\city{Cambridge}, \postcode{CB2 1PZ}, \country{United Kingdom}}}


\abstract{Shape optimisation of thin-shell structures requires a flexible, differentiable geometric representation suitable for gradient-based optimisation.  We propose a neural parametric geometry representation (NRep) for shells based on a neural network with periodic activation functions. The NRep is defined using a multi-layer perceptron (MLP), which maps the parametric coordinates of mid-surface vertices to their physical coordinates. A structural compliance optimisation problem is posed to optimise the shape of a thin-shell parameterised by the NRep subject to a volume constraint, with the network parameters as design variables. The resulting shape optimisation problem is solved using a gradient-based optimisation algorithm. Benchmark examples with classical solutions and comparisons with the free-form deformation method demonstrate that the proposed NRep is capable of representing shell geometries with local geometric features using a small set of network parameters. The robustness of the approach has been demonstrated with different initial geometries, boundary conditions and neural network hyperparameters. The approach also exhibits potential for complex lattice-skin structures, owing to the compact and expressive geometry representation afforded by the NRep.}

\keywords{thin shells, shape optimisation, neural network, neural parametric representation, lattice-skin structures}

\maketitle

\section{Introduction} \label{sec:intro}

Thin-shell structures are widely used in engineering because of their high stiffness-to-weight ratio and superior load-carrying capacity. With a thickness much smaller than the mid-surface dimensions, the structural response of a thin-shell is highly sensitive to its mid-surface geometry. Small geometric changes to the mid-surface can result in large and often non-local variations in structural behaviour. Conventional thin-shell shape optimisation approaches, updating the finite element mesh that approximates the mid-surface, are problematic for practical applications. In problems involving significant shape changes, the mesh may become distorted or entangled, and this in turn adversely affects both finite element analysis and sensitivity analysis required by gradient-based optimisation algorithms. Moreover, the geometries of the mid-surface and the finite element model are coupled but usually represented in different forms, so converting an optimised mesh result back to the original geometric representation introduces additional approximation errors. In the context of isogeometric analysis, where geometry and finite element discretisation use the same spline basis~\cite{hughes2005isogeometric, cirak2002integrated}, shape updates are intrinsically linked to the analysis representation. This motivates the development of optimisation strategies based on alternative surface parameterisations.

Parametric representations are the standard in CAD systems for geometric modelling, including meshes and splines (for example, B-splines, NURBS and subdivision surfaces). Shape optimisation methods built on these parameterisation techniques therefore predominate. Design variables are commonly defined as the positions of mesh vertices~\cite{noboru1986adaptive} or control points of spline surfaces~\cite{braibant1984shape}. However, direct optimisation of vertex positions in a mesh is prone to generate unrealistic designs with jagged geometries, as moving an individual vertex provides only local control over the mesh and may cause the optimisation to become trapped in local minima. To mitigate these problems, boundary regularity constraints, smooth remeshing and filtering techniques have been proposed to enforce smoothness in mesh-based shape optimisation~\cite{le2011gradient, hojjat2014vertex, swartz2023yet}. Spline surfaces are less prone to such geometric artefacts due to their higher-order spline basis functions and broader local supports, but undesirable geometries may still arise when control points are adjusted. Several techniques have been proposed to improve shape optimisation results of specific spline surfaces, for example, augmenting the design space by optimising both the positions and weights of control points of NURBS surfaces~\cite{qian2010full}, and adopting multiresolution schemes for subdivision surfaces to regularise the shape optimisation to prevent sub-optimal jagged geometries~\cite{bandara2015boundary, bandara18Isogeometric}. Alternatively, the vertex-based design variables can be linked indirectly to geometries using global interpolation techniques, such as radial basis functions~\cite{jakobsson2007mesh, rendall2008unified} and free-form deformation~\cite{samareh2004aerodynamic, xiao2022infill, xiao2022free, zhao2024automated}, enabling shape control independent of the underlying surface parameterisation.

Implicit representations provide a fundamentally different paradigm for geometric modelling. An implicit surface is usually defined as the zero level set of a scalar field, which offers inherent smoothness and flexibility with respect to topological changes. Level set methods based on implicit representations have been widely applied in structural optimisation~\cite{sethian2000structural, wang2003level, kobayashi2024shell}, where the level set is discretised on a regular grid and structural boundaries evolve according to shape sensitivities or topological derivatives~\cite{allaire2004structural}. More recently, neural distance fields based on deep learning have been introduced~\cite{park2019deepsdf, chibane2020neural, sitzmann20Implicit}, in which a continuous distance field is encoded by a neural network and its zero level defines an implicit surface, commonly referred to as a neural implicit surface. Signed distance fields require inside/outside labelling and are therefore restricted to orientable and closed surfaces. Surface meshes can subsequently be extracted from the fields using established algorithms, such as Marching Cubes or Dual Contouring. For open surfaces with boundaries, unsigned distance fields (UDFs) are required; however, accurate surface mesh extraction from UDFs is more challenging, as the absence of sign changes and non-negligible approximation errors near the surface, particularly at sharp edges, can lead to reduced fidelity~\cite{zhang2023surface}. Neural implicit representations have also been incorporated into topology optimisation, where spatial coordinates are mapped to a continuous distance field~\cite{zhang23Topology, hu2024if}, analogous to the discrete level sets in conventional level set methods. Alternatively, when topology optimisation is formulated in terms of physical densities, a neural network can directly generate a continuous spatial density field~\cite{hoyer2019neural, chandrasekhar21Topology, zehnder2021ntopo}, resembling the element densities in the SIMP method~\cite{bendsoe1999material}. In such neural optimisation frameworks, the optimisation problem is reparameterised in terms of neural network parameters, and this could lead to improved optimisation results~\cite{sanu2025neural} given appropriate network architectures and training strategies.

Nonetheless, using implicit surfaces, whether defined with grid-based level sets or continuous neural fields, cannot precisely represent thin-shells with open boundaries and sharp features that are commonly present in engineering applications. In contrast, parametric surfaces can naturally represent these geometries, which explains why thin-shell shape optimisation using parametric surfaces has been explored extensively over the past decades. However, learning-based approaches for parametric surfaces remain scarce compared with implicit representations~\cite{morreale2021neural}, and even more so for parametric shape optimisation. In this work, we propose a framework for thin-shell shape optimisation based on neural parametric representations (NReps). The proposed pipeline integrates NReps, finite element analysis, nonconvex optimisation and automatic differentiation in a unified manner. We first parameterise the mid-surface of the thin-shell by a neural network, specifically a multi-layer perceptron (MLP), which maps global parametric coordinates to physical coordinates. The thin-shell shape is therefore parameterised by the network parameters (weights and biases). We use sinusoidal activation functions in the hidden layers to improve the representation of high-frequency geometric details~\cite{sitzmann20Implicit}. Given the parameterised mid-surface, a standard thin-shell finite element analysis is performed to compute the displacement field. The constrained optimisation problem is then formulated in the space of neural network parameters and solved with a gradient-based optimiser, such as sequential quadratic programming (SQP). Analytical gradients are obtained by chaining finite element sensitivities with automatic differentiation of the neural network. Updating the network parameters directly modifies the mid-surface shape, and the optimisation process is driven by the objective function, which is structural compliance minimisation in the presented examples.

Compared with spline-based representations, which rely on fixed piecewise polynomials defined by knot vectors and control meshes, NReps are essentially meshless and capable of approximating highly general parametric mappings. Owing to the universal approximation property of neural networks~\cite{Nielsen2015neural, Goodfellow2016deep}, NReps can represent non-polynomial geometries that would otherwise require many spline pieces. Furthermore, the nonlinearity of NReps is easy to tune through the neural network architecture, including the depth and width of hidden layers and the choice of activation functions, giving direct control over approximation accuracy and smoothness. In particular, sinusoidal activations provide well-defined gradients as well as higher-order derivatives~\cite{sitzmann20Implicit}, which benefits both network training and gradient-based shape optimisation. In the proposed framework, the number of design variables in the NRep is independent of the finite element discretisation. Consequently, changing the NRep that parameterises the same mid-surface does not require the finite element model to be remeshed or modified, which distinguishes the proposed framework from conventional spline or mesh based approaches.

The remainder of the paper is organised as follows. In Section 2, the thin-shell finite element analysis based on mid-surface geometry is briefly introduced. Section 3 describes the neural parametric representations. Building upon neural parametric representations, we introduce the shape optimisation framework and derive the sensitivity analysis in Section 4. Section 5 presents optimisation examples with several thin-shell shapes, topologies and boundary conditions, demonstrating the effectiveness and robustness of the proposed approach.

\section{Review of thin-shell finite element analysis} \label{sec:shellFea}
The Kirchhoff--Love thin-shell is considered, with the mid-surface $\Omega$ parameterised by curvilinear coordinates $\vec{\eta} = (\eta^1, \eta^2) \in \mathbb{R}^2$. The physical coordinates of points on the mid-surface are denoted by $\vec{x}(\vec{\eta}) \in \mathbb{R}^3$. The covariant and contravariant basis vectors of the mid-surface are defined as
\begin{equation}
	\vec{a}_\alpha = \frac{\partial \vec{x}}{\partial \eta^\alpha}
	\quad \text{and} \quad
	\vec{a}^\beta \cdot \vec{a}_\alpha = \delta_\alpha^\beta \, ,
\end{equation}
respectively, where $\delta_\alpha^\beta$ is the Kronecker delta, and Greek indices take values $\{1, 2\}$. The unit normal $\vec{a}_3$ is given by
\begin{equation}
	\vec{a}_3 = \frac{\vec{a}_1 \times \vec{a}_2}{|\vec{a}_1 \times \vec{a}_2|} \, .
\end{equation}
For a deformed thin-shell with mid-surface displacement $\vec{u}$, the total potential energy under an external loading $\vec{f}$ is
\begin{equation}
\begin{split}
	\Pi(\vec{u}) &= \int_{\Omega} \left(W^{\text{m}}(\vec{\alpha}(\vec{u})) + W^{\text{b}}(\vec{\beta}(\vec{u}))\right) \mathrm{d}\Omega \\
	&\quad - \int_\Omega \vec{f}\cdot\vec{u}\,\mathrm{d}\Omega \, ,
\end{split}
\end{equation}
where $W^{\text{m}}$ and $W^{\text{b}}$ denote the membrane and bending strain energy densities, respectively. The membrane strain $\vec{\alpha}$ and the bending strain $\vec{\beta}$ can be expressed in terms of the covariant and contravariant basis vectors (refer to the monograph~\cite{ciarlet05Introduction} for more details). In the case of small displacements, the linearised strains are given by~\cite{cirak00Subdivision}
\begin{equation}
	\vec{\alpha} = \frac{1}{2}\left(\vec{a}_{\alpha} \cdot \vec{u}_{,\beta} + \vec{u}_{,\alpha} \cdot \vec{a}_{\beta}\right)\vec{a}^{\alpha}\otimes\vec{a}^{\beta}
\end{equation}
and
\begin{equation}
\begin{split}
	\vec{\beta} &= \Big(-\vec{u}_{,\alpha\beta} \cdot \vec{a}_3 + \frac{1}{\sqrt{a}}\big[\vec{u}_{,1} \cdot (\vec{a}_{\alpha,\beta} \times \vec{a}_2) \\
	&+ \vec{u}_{,2} \cdot (\vec{a}_1 \times \vec{a}_{\alpha,\beta})\big] + \frac{\vec{a}_3 \cdot \vec{a}_{\alpha,\beta}}{\sqrt{a}}\big[\vec{u}_{,1} \cdot (\vec{a}_2 \times \vec{a}_3) \\
	&+ \vec{u}_{,2} \cdot (\vec{a}_3 \times \vec{a}_1)\big]\Big) \vec{a}^{\alpha} \otimes \vec{a}^{\beta} \, ,
\end{split}
\end{equation}
where $\sqrt{a} = |\vec{a}_1 \times \vec{a}_2|$ and the commas denote partial derivatives with respect to the curvilinear coordinates.

The variational equilibrium equations are obtained from the stationary points of the potential energy as
\begin{equation} \label{eq:stationary}
\begin{split}
	\frac{\partial \Pi(\vec{u})}{\partial \vec{u}} \delta \vec{u} &= \int_\Omega\left(\frac{\partial W^\text{m}(\vec{\alpha})}{\partial \vec{u}}\delta\vec{u} + \frac{\partial W^\text{b}(\vec{\beta})}{\partial \vec{u}}\delta\vec{u}\right)\mathrm{d}\Omega \\
	&\quad - \int_\Omega \vec{f}\delta\vec{u}\,\mathrm{d}\Omega = \vec{0} \, .
\end{split}
\end{equation}
The position vector $\vec{x}$ and the displacement field $\vec{u}$ of the mid-surface are discretised using quantities defined at discrete vertices, i.e.,
\begin{equation}
	\vec{x} \approx \sum_i B_i(\vec{\eta})\mat{x}_i \, , \quad
	\vec{u} \approx \sum_i B_i(\vec{\eta})\mat{u}_i \, ,
\end{equation}
where $B_i(\vec{\eta})$ are basis functions (e.g. Catmull--Clark subdivision or other spline basis functions) evaluated at $\vec{\eta}$, and $\mat{x}_i$ and $\mat{u}_i$ are vertex coordinates and displacements, respectively. Following discretisation of the variational formulation~\eqref{eq:stationary}, the discretised equilibrium equation system is obtained:
\begin{equation}
	\mat{K}\mat{u} = \mat{f}
\end{equation}
with the discretised global stiffness matrix $\mat{K}$, displacement vector $\mat{u}$ and force vector $\mat{f}$.

\section{Neural parametric representation (NRep)} \label{sec:neuRep}
For a surface $\Omega$, consider a neural network $\mathcal{N}_{\vec{\theta}}(\vec{\overline{\vec{x}}})$ that takes as input the parametric coordinates $\overline{\vec{x}} \in \mathbb{R}^2$ of points, with $\vec{\theta}$ denoting the neural network parameters. The output of the network is the physical coordinates of the points $\vec{x} \in \mathbb{R}^3$, i.e.,
\begin{equation} \label{eq:nrep}
	\vec{x} = \mathcal{N}_{\vec{\theta}}(\overline{\vec{x}}) \, ,
\end{equation}
which defines a neural parametric representation of the surface. The surface is represented by the neural network, rather than an explicit parametric form, such as a triangulation or a spline surface.

\begin{figure}
	\centering
	\includegraphics[width=\linewidth]{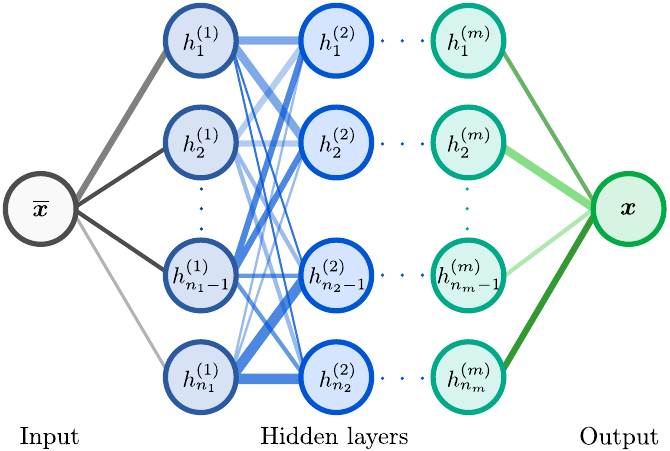}
	\caption{Multi-layer perceptron for neural parametric representation.}
	\label{fig:neuralNetwork}
\end{figure}

A simple multi-layer perceptron (MLP) is adopted as the neural network architecture, as illustrated in Figure~\ref{fig:neuralNetwork}. The network consists of an input layer, an output layer and a few hidden layers. The input layer takes the parametric coordinates of points on the mid-surface, and the output layer gives their coordinates in the physical space. This establishes a mapping from the parametric space to the physical space. The parametric coordinates are chosen such that $\overline{\vec{x}} \in [0, 1]\times[0, 1]$. It is worth noting that the input of the neural network is a single point in the parametric space and the output is the corresponding point on the surface. The same set of network parameters defines a continuous mapping and is therefore used to represent all points on the surface. For a discretised mid-surface with $n_v$ vertices, the network is evaluated simultaneously at all vertices for computational efficiency.

The network parameters $\vec{\theta}$ comprise the weights $\vec{w}$ and biases $\vec{b}$ of the neural network, which determine the activations of nodes in the hidden layers and are determined through network training via backpropagation. The activations of the $i$-th hidden layer are computed as
\begin{equation}
	\vec{c}_i = \vec{W}_i \vec{h}^{(i - 1)} + \vec{b}_i \, ,
\end{equation}
where $\vec{W}_i$ and $\vec{b}_i$ are the weight matrix and bias vector of the $i$-th hidden layer, respectively; $\vec{h}^{(i-1)}$ denotes transformed quantities of the previous hidden layer via an activation function. To enhance the representation of fine geometric details, rather than using conventional ReLU or Tanh activation functions alone, we adopt the activation function in periodic forms~\cite{sitzmann20Implicit, zhang23Topology} as follows,
\begin{equation} \label{eq:periodicActivation}
	\vec{h}^{(i)}(\vec{c}_i) = \sin(\omega \vec{c}_i + \vec{\delta})
\end{equation}
where $\omega$ is a period-related coefficient and $\vec{\delta}$ is a phase shift, both of which are prescribed and untrainable in the study. The sinusoidal activation function is applied componentwise in~\eqref{eq:periodicActivation}.

\begin{figure}
	\centering
	\begin{subfigure}{\linewidth}
		\centering
		\includegraphics[width=0.8\linewidth]{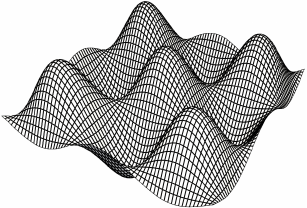}
		\caption{Surface mesh defined with $5\cos(x/2)\cos(y/2)$.}
		\label{fig:surfMesh}
	\end{subfigure}
	\begin{subfigure}{\linewidth}
		\centering
		\includegraphics[width=0.45\linewidth]{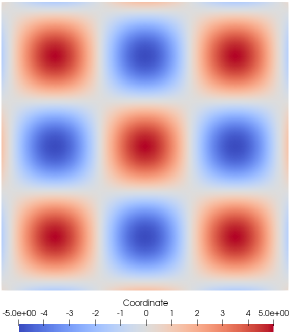}
		\hfill
		\includegraphics[width=0.45\linewidth]{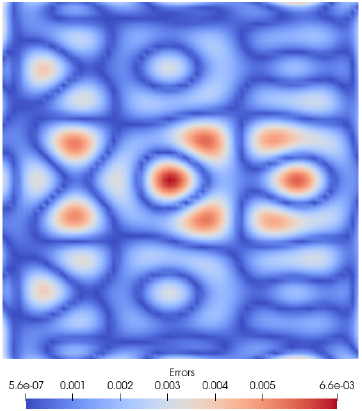}
		\caption{Z-coordinates obtained from NRep (left) and relative errors of vertex coordinates (right), with a range between $5.6 \times 10^{-7}$ and $6.6 \times 10^{-3}$.}
		\label{fig:approxError}
	\end{subfigure}
	\begin{subfigure}{\linewidth}
		\centering
		\includegraphics[width=\linewidth]{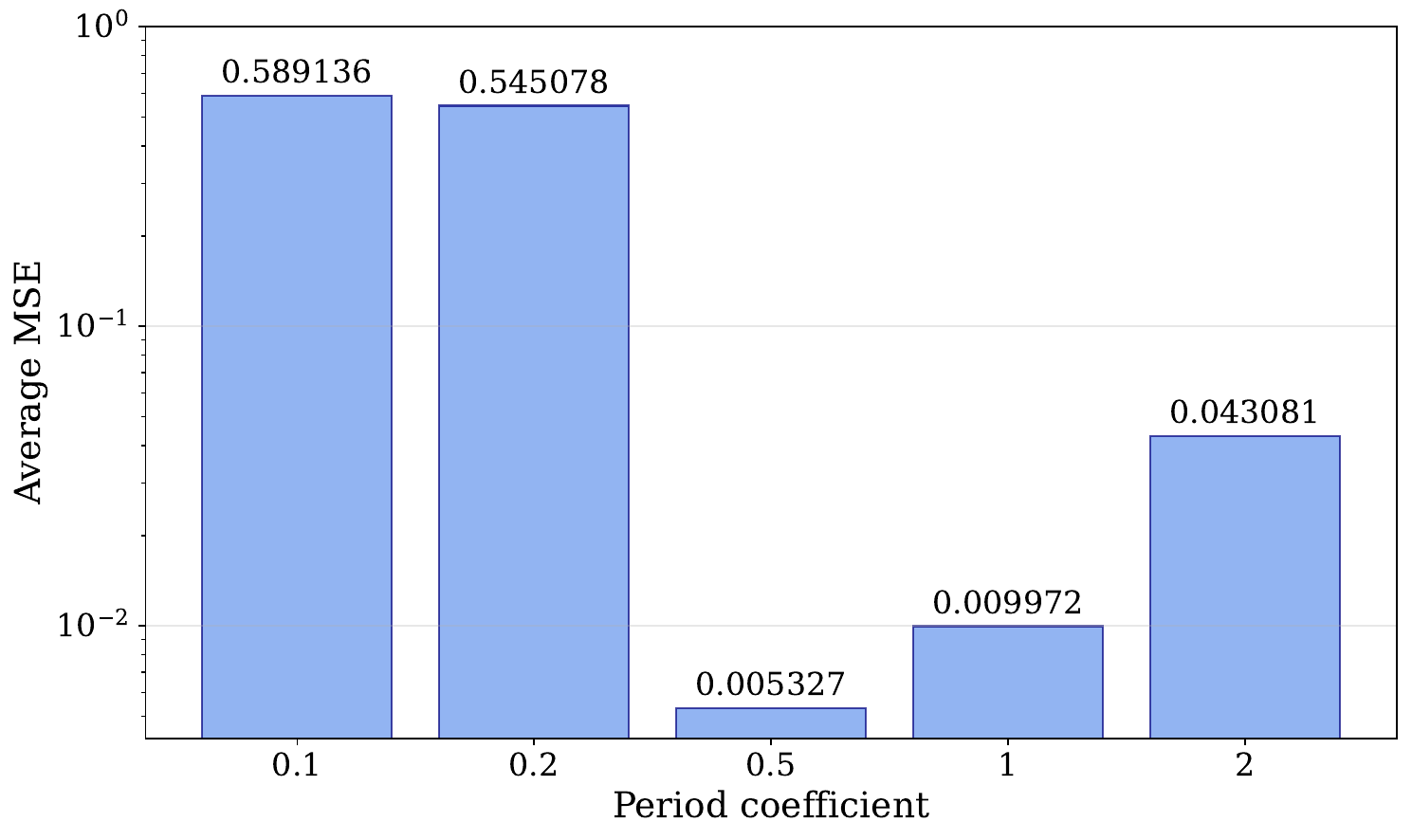}
		\caption{Influence of period coefficient $\omega$ on approximation error of NRep.}
		\label{fig:periodEffect}
	\end{subfigure}
	\caption{Surface approximation with NRep.}
	\label{fig:NRepApprox}
\end{figure}

Figure~\ref{fig:NRepApprox} shows an example of approximating a free-form surface with NRep. The projected domain of the surface is a square of size $20 \times 20$, and the corresponding surface mesh consists of $4096$ elements and $4225$ vertices. The neural network of NRep has three hidden layers, each containing five nodes, with activation functions of the form $\sin(\vec{c}_i/2 + \pi/4)$, resulting in a total of $261$ trainable parameters. The network input is the projected coordinates of the surface scaled to the range $[0, 1]\times[0, 1]$. The NRep of the surface is obtained by minimising the mean squared error (MSE) between the approximated and reference coordinates of the mesh vertices, using the Adam optimiser with a learning rate of $0.01$. After $20{,}000$ training epochs, the MSE reaches $1.277 \times 10^{-3}$. The relative errors of the vertex coordinates are shown in Figure~\ref{fig:approxError}, with a maximum error of $0.66\%$ relative to the surface dimension, demonstrating the high approximation accuracy of NRep. The influence of the period coefficient $\omega$ in the activation function is further investigated by considering different values. For each value of $\omega$, the NRep is trained independently ten times with the same $20{,}000$ epochs, and the average MSE is reported. The comparison results, shown in Figure~\ref{fig:periodEffect}, indicate that $\omega = 0.5$, which closely matches the frequency of the original surface, corresponds to an NRep with the smallest approximation error. Smaller values of $\omega$ result in higher approximation errors when representing surfaces with fine geometric details.

\section{Thin-shell shape optimisation} \label{sec:shapeOpt}
\subsection{Optimisation scheme}
\begin{figure*}
	\centering
	\includegraphics[width=0.9\linewidth]{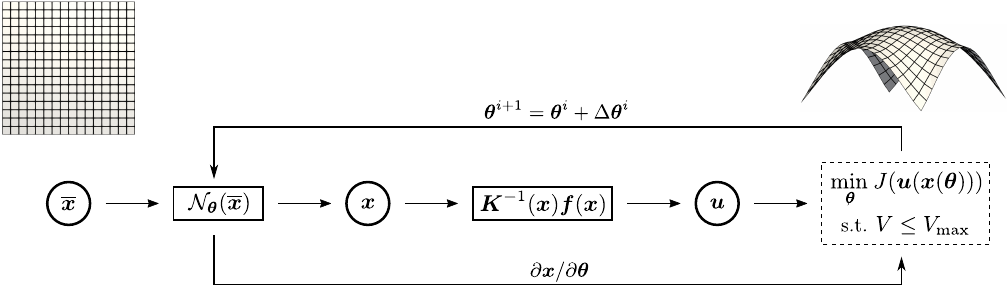}
	\caption{Flowchart of thin-shell shape optimisation with NRep.}
	\label{fig:shapeOptScheme}
\end{figure*}

In shape optimisation we aim to find the optimal shape of the thin-shell, which minimises the structural compliance while satisfying a prescribed volume constraint. Since the mid-surface of the thin-shell is parameterised by the neural network, the network parameters $\vec{\theta}$ are the design variables. The optimised thin-shell shape is inferred readily via~\eqref{eq:nrep} after the network parameters are determined. The shape optimisation problem can be formulated as follows,
\begin{subequations} \label{eq:shapeProb}
\begin{align}
	\min_{\vec{\theta}} \quad &J(\vec{x}(\vec{\theta})) = \mat{u}^\trans(\vec{x})\mat{K}(\vec{x})\mat{u}(\vec{x}) \, , \\
	\text{s.t.} \quad &\mat{K}(\vec{x})\mat{u}(\vec{x}) = \mat{f} \, , \label{eq:staticEquilibrium} \\
	& V(\vec{x}) \leq V_{\max} \, , \label{eq:volCons}
\end{align}
\end{subequations}
where $V$ and $V_{\max}$ are the thin-shell volume and the prescribed maximum volume, respectively.

One approach to optimising the network parameters is to train the neural network using the optimisers available in neural network frameworks (e.g. the Adam optimiser), which primarily deal with unconstrained optimisation problems. In such cases, the loss function of the neural network is penalised to incorporate the volume constraint~\eqref{eq:volCons}~\cite{chandrasekhar21Topology, zhang23Topology}, such that the original shape optimisation~\eqref{eq:shapeProb} is converted into an unconstrained problem by minimising the loss function directly. However, this approach is sensitive to the choice of penalty parameter, and the available optimisers are generally not well suited for shape optimisation.

In our implementation, the network parameters are not optimised by stochastic gradient descent as in the conventional neural network training. Instead, they are handled within a separate shape optimisation module and updated iteratively, as illustrated in Figure~\ref{fig:shapeOptScheme}. Gradient-based optimisation algorithms, such as sequential quadratic programming (SQP) and the method of moving asymptotes (MMA), can be used for solving the constrained shape optimisation problem~\eqref{eq:shapeProb}.

Before shape optimisation, the network parameters are initialised by fitting the NRep to the initial mid-surface geometry. This pre-training stage minimises the MSE loss function with the Adam optimiser, as described in Figure~\ref{fig:NRepApprox}. The network initialisation is performed separately from the subsequent shape optimisation process, although both stages use the same neural network architecture. The pre-trained network parameters are used as the starting point for the subsequent shape optimisation. The purpose of the pre-training stage is to provide a geometric representation that is consistent with the initial thin-shell shape and to avoid starting the shape optimisation from randomly initialised network parameters, which may adversely affect the optimisation process. It should be noted that the pre-training process itself starts from randomly initialised parameters. Consequently, different initialisations may lead to different sets of network parameters, even though they all approximate the same initial geometry.

At each optimisation iteration during shape optimisation, the sensitivities with respect to the network parameters are evaluated, and the parameters are updated according to the chosen optimisation algorithm. The process is repeated until the optimisation termination criterion is satisfied, after which the optimised thin-shell shape is inferred from the neural network output.

\subsection{Sensitivity analysis for NRep update}
During optimisation the derivatives of the objective and constraint functions with respect to the network parameters $\vec{\theta}$ are required. Assuming a constant external force vector $\mat{f}$, the derivative of the objective function can be expressed as
\begin{equation}
	\frac{\partial J}{\partial \theta_j} = \sum_i \frac{\partial J}{\partial x_i} \frac{\partial x_i}{\partial \theta_j} = \sum_i \left(\mat{u}^\trans \frac{\partial \mat{K}}{\partial x_i} \mat{u} + 2 \mat{f}^\trans \frac{\partial \mat{u}}{\partial x_i} \right) \frac{\partial x_i}{\partial \theta_j}
\end{equation}
and considering the static equilibrium~\eqref{eq:staticEquilibrium} yields
\begin{equation}
	\frac{\partial J}{\partial \theta_j} = -\sum_i \left(\sum_e \mat{u}_e^\trans \frac{\partial \mat{K}_e}{\partial x_i} \mat{u}_e \right) \frac{\partial x_i}{\partial \theta_j} \, ,
\end{equation}
where $\mat{K}_e$ and $\mat{u}_e$ are the element stiffness matrix and the element displacement vector, respectively; for Kirchhoff--Love thin-shells, the derivatives $\partial \mat{K}_e / \partial x_i$ can be computed directly~\cite{bandara18Isogeometric}; the derivatives $\partial x_i / \partial \theta_j$ are obtained by automatic differentiation of the neural network.

Similarly, the sensitivity of the thin-shell volume constraint is given by
\begin{equation}
	\frac{\partial V}{\partial \theta_j} = \sum_i t\frac{\partial A}{\partial x_i} \frac{\partial x_i}{\partial \theta_j} \, ,
\end{equation}
where $A$ denotes the mid-surface area and $t$ is the shell thickness.

\section{Examples} \label{sec:examples}
The application and effectiveness of the proposed shape optimisation approach are demonstrated with four numerical examples. The first two examples consider strip and square thin-shells with available benchmark optimisation results. The third example investigates shape optimisation of thin-shells with varying boundary conditions and topologies. Finally, a complex lattice-skin geometry is generated based on an optimised NRep of a thin-shell in the last example. In all examples the thin-shells are optimised only in the height direction relative to their planform. Consequently, the output of the NRep is a scalar and its input is a two-dimensional vector. The neural networks are implemented in TensorFlow~2 and the SQP optimisation algorithm provided in the NLopt library is adopted~\cite{Steven2020nlopt}.

\subsection{Catenary strip}
First, a thin-shell strip with two short edges pinned is considered. The projected area of the strip is $20 \times 1$, and the thickness is $0.1$. A uniform load of magnitude $10$ is applied on the strip. The Young's modulus and the Poisson's ratio of the material are set to $7 \times 10^7$ and $0.35$, respectively. A volume constraint is imposed such that the optimised strip has a volume $5\%$ greater than the initial volume.

The initial strip is discretised with $99$ vertices and $64$ quadrilateral shell elements. The neural network consists of two hidden layers with five nodes each. The number of network parameters is $51$. Three sets of activation functions are considered for the two hidden layers: both with ReLU functions (ReLU + ReLU), one with a periodic function and the other with ReLU (periodic + ReLU), and both with periodic functions (periodic + periodic). The periodic activation function takes the form $\sin(\vec{c}_i + \pi/4)$, cf.~\eqref{eq:periodicActivation}.

\begin{figure*}
	\centering
	\begin{subfigure}{\linewidth}
		\includegraphics[width=0.32\linewidth]{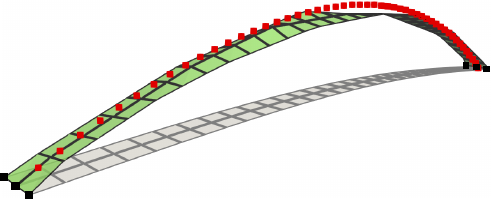}
		\hfill
		\includegraphics[width=0.32\linewidth]{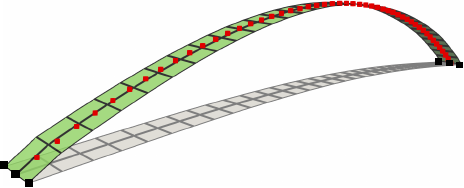}
		\hfill
		\includegraphics[width=0.32\linewidth]{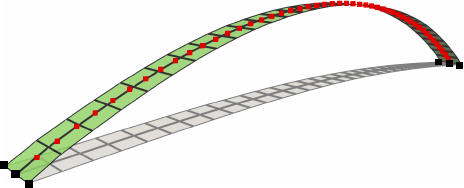}
		\caption{Optimised shapes with red points discretising the catenary curve. From left to right, the mean squared errors relative to the catenary curve are $2.82 \times 10^{-2}$, $1.35 \times 10^{-3}$ and $5.21 \times 10^{-5}$, respectively.}
	\end{subfigure}
	\begin{subfigure}{\linewidth}
		\includegraphics[width=0.32\linewidth]{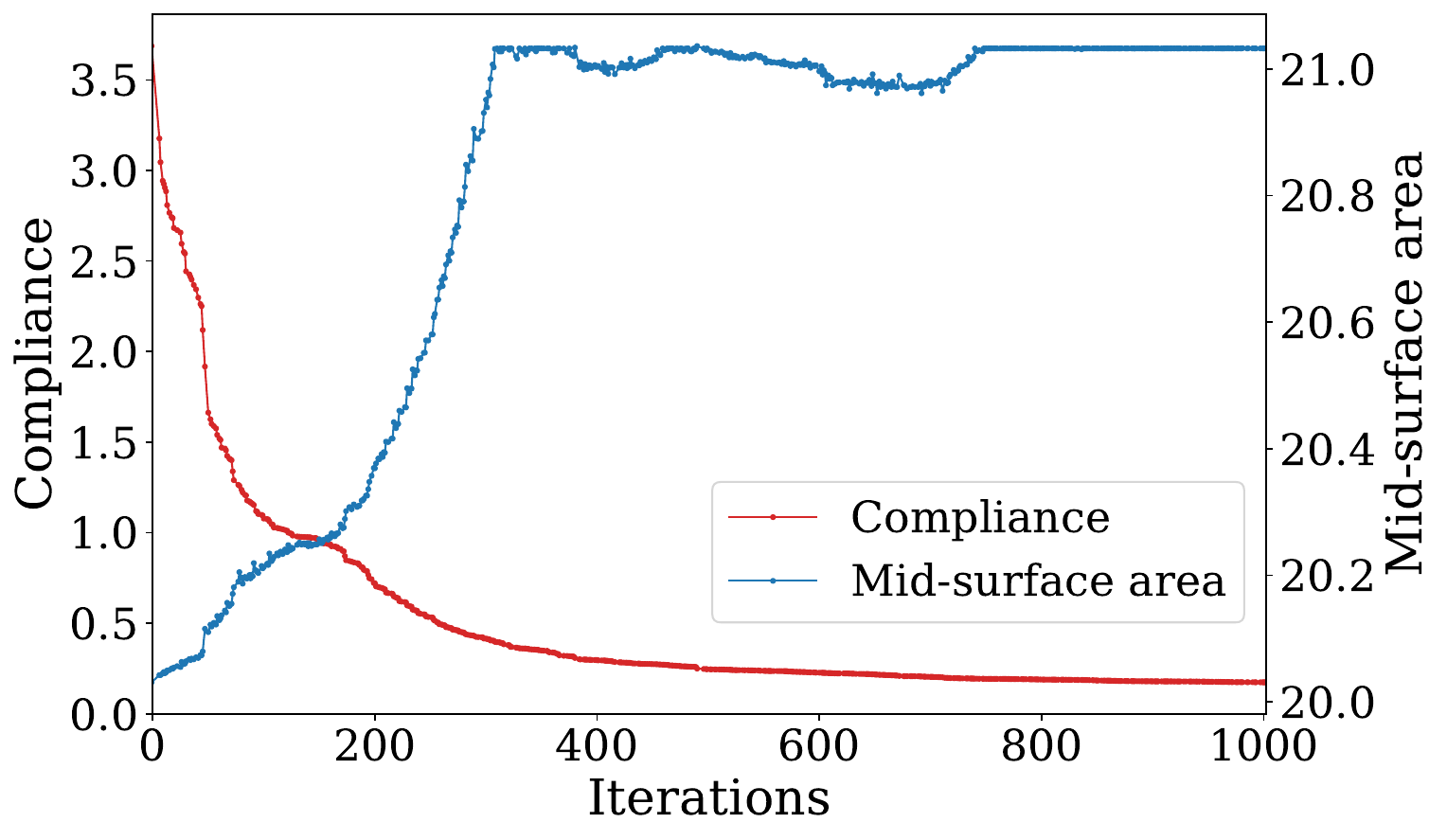}
		\hfill
		\includegraphics[width=0.32\linewidth]{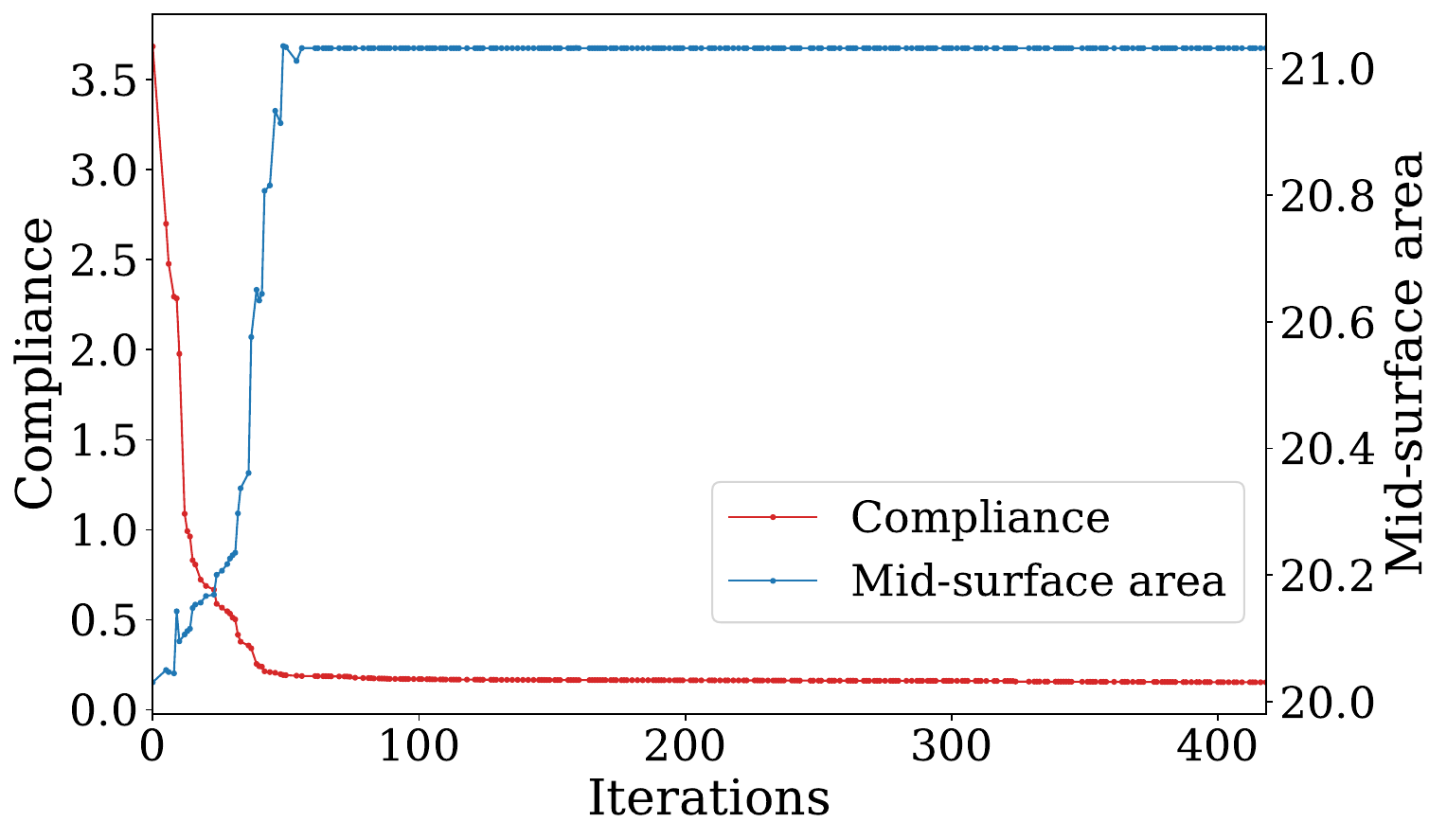}
		\hfill
		\includegraphics[width=0.32\linewidth]{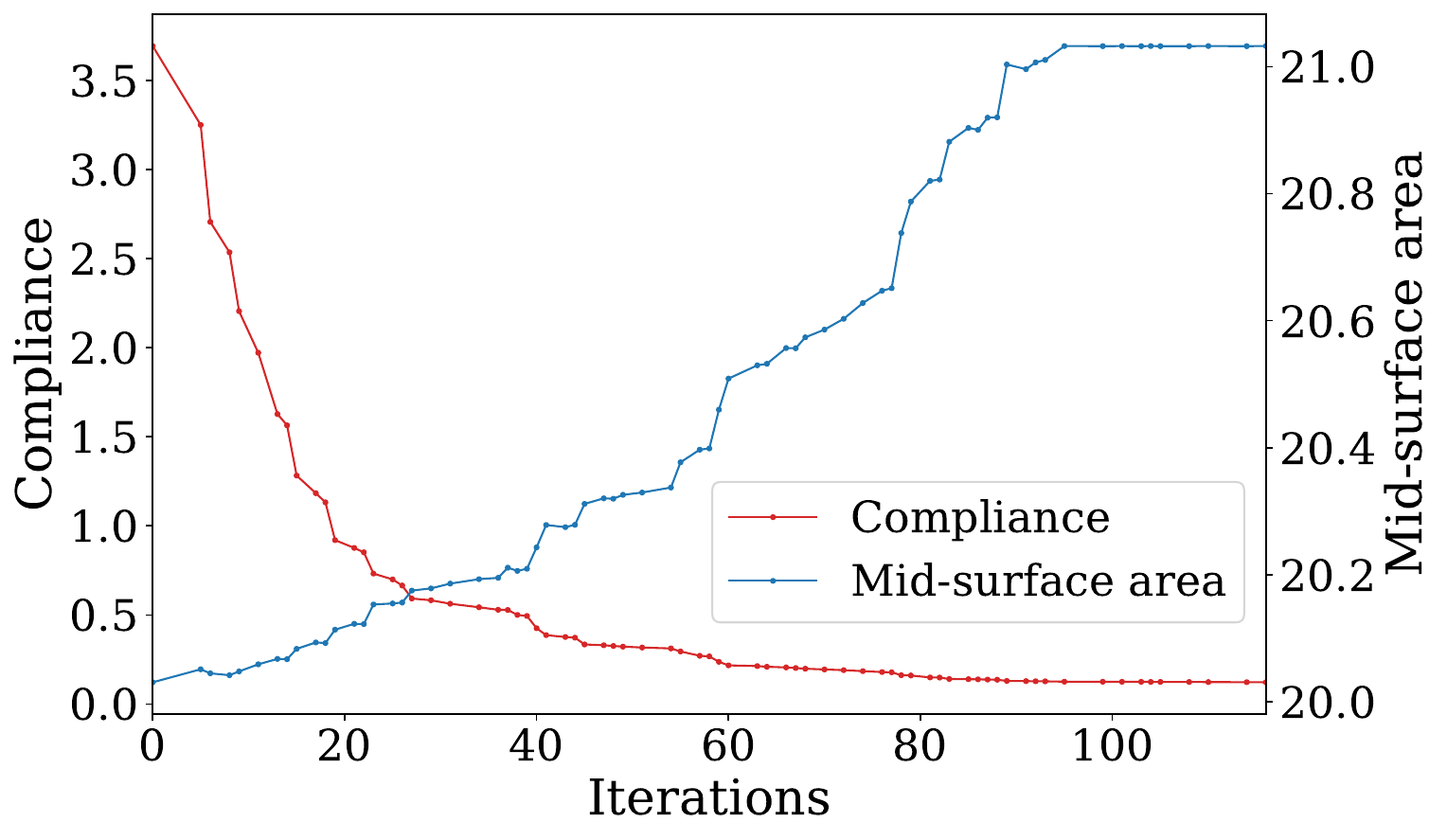}
		\caption{Compliance and mid-surface area variations during optimisation. From left to right, the optimised compliance values are $0.164$, $0.143$ and $0.113$, respectively.}
	\end{subfigure}
	\caption{Shape optimisation of the strip under uniform loading, using ReLU + ReLU (left), periodic + ReLU (middle) and periodic + periodic (right) activation functions. For data visualisation, the iterations in diagrams (b) are truncated to the compliance values that are no more than 0.01 larger than the optimised compliance.}
	\label{fig:archShapeOpt}
\end{figure*}

The resulting optimised strip shapes obtained with different activation functions are shown in Figure~\ref{fig:archShapeOpt} and are compared with the analytical catenary curve, which represents the optimal shape of a curve under self-weight~\cite{gohnert23Catenary}. The comparison indicates that using periodic activation functions yields shapes closer to the catenary curve, while requiring fewer optimisation iterations.

To investigate the influence of the number of network parameters on the optimisation result, we consider hidden layers with two or ten nodes, corresponding to $15$ and $151$ network parameters, respectively. Both architectures use the ReLU activation function in their two hidden layers. The optimisation results obtained with the three neural networks are summarised in Table~\ref{tab:stripCellCompare}. With more network parameters, the optimised shape changes towards the catenary curve with smaller compliance, due to the increased flexibility of shape representation, albeit at the cost of additional optimisation iterations.

\begin{table}
	\centering
	\caption{Influence of the number of nodes in each hidden layer on shape optimisation using only ReLU activation functions.}
	\begin{tabular}{cccc}
		\toprule
		\# nodes & \# parameters & MSE & compliance \\
		\midrule
		$2$ & $15$ & $3.88 \times 10^{-2}$ & 0.711 \\
		$5$ & $51$ & $2.82 \times 10^{-2}$ & 0.164 \\
		$10$ & $151$ & $1.48 \times 10^{-3}$ & 0.123 \\
		\bottomrule
	\end{tabular}
	\label{tab:stripCellCompare}
\end{table}

\subsection{Thin-shell roof}
Next, a square thin-shell under a uniform load of magnitude $10$ is considered for shape optimisation. The four mid-edge points of the thin-shell are fixed. The projected area of the thin-shell is $20 \times 20$, and the thickness is $0.1$. The Young's modulus and the Poisson's ratio of the material are set to $7 \times 10^7$ and $0.35$, respectively. The volume is constrained to be no greater than $1.2$ times the initial volume.

\begin{figure*}
	\centering
	\begin{subfigure}{\linewidth}
		\includegraphics[width=0.32\linewidth]{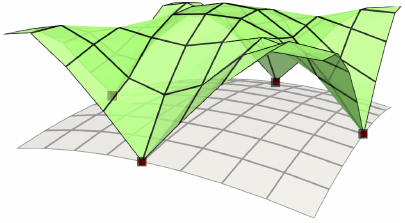}
		\hfill
		\includegraphics[width=0.32\linewidth]{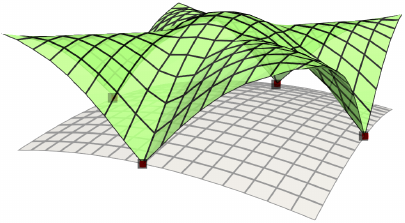}
		\hfill
		\includegraphics[width=0.32\linewidth]{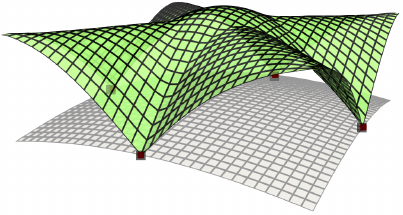}
		\caption{Optimised shapes with $64$ (left), $256$ (middle) and $1024$ (right) shell elements.}
	\end{subfigure}
	\\
	\begin{subfigure}{\linewidth}
		\includegraphics[width=0.32\linewidth]{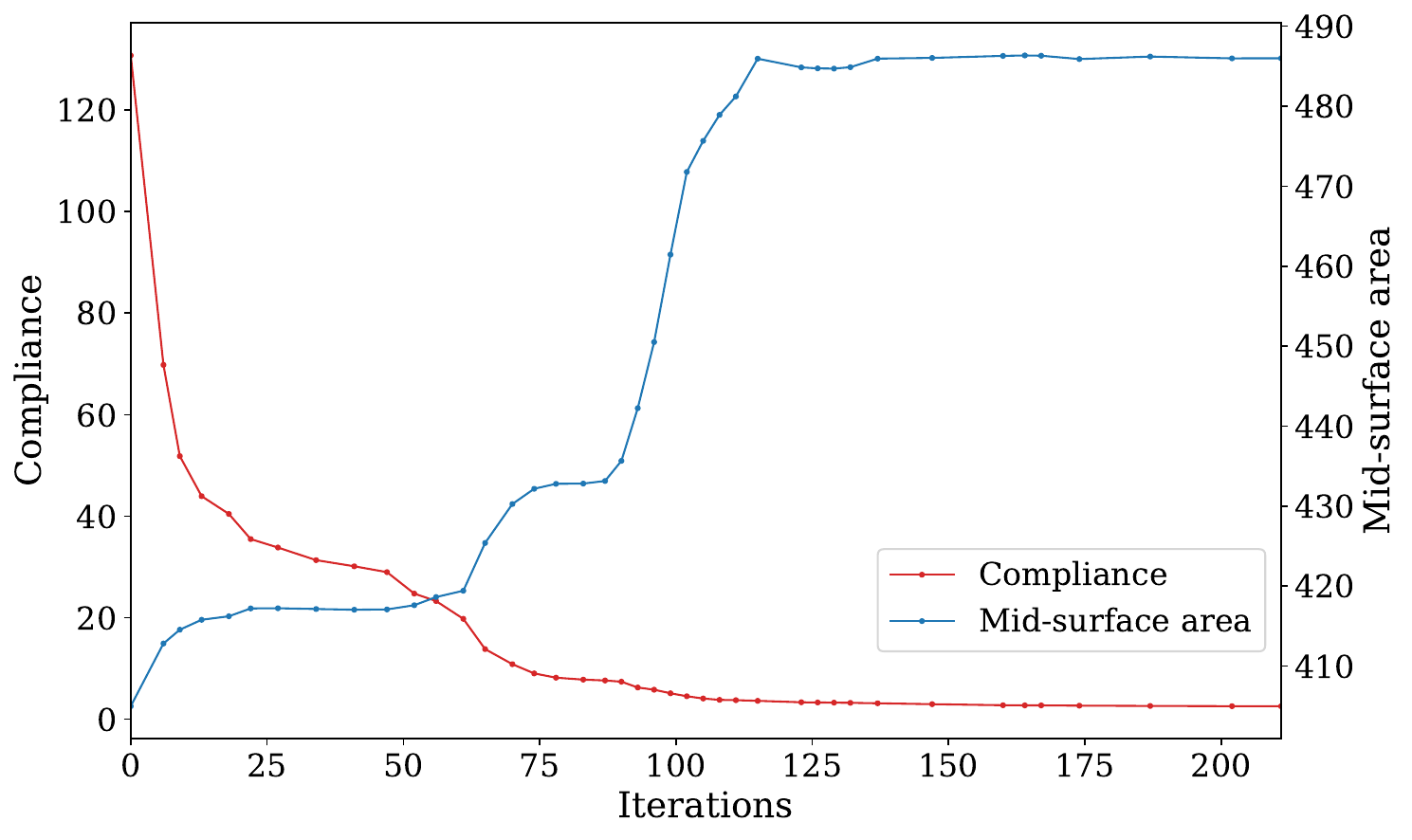}
		\hfill
		\includegraphics[width=0.32\linewidth]{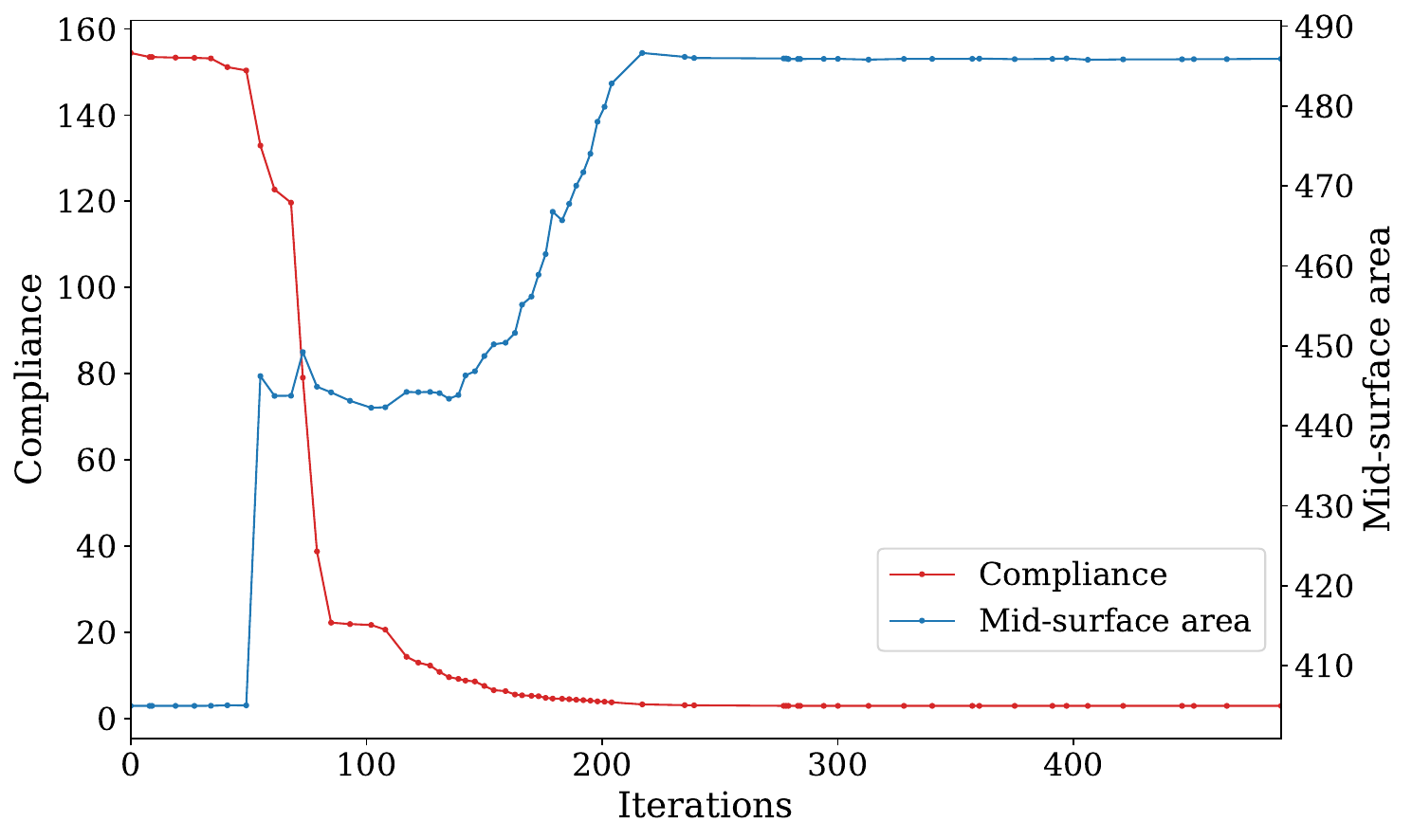}
		\hfill
		\includegraphics[width=0.32\linewidth]{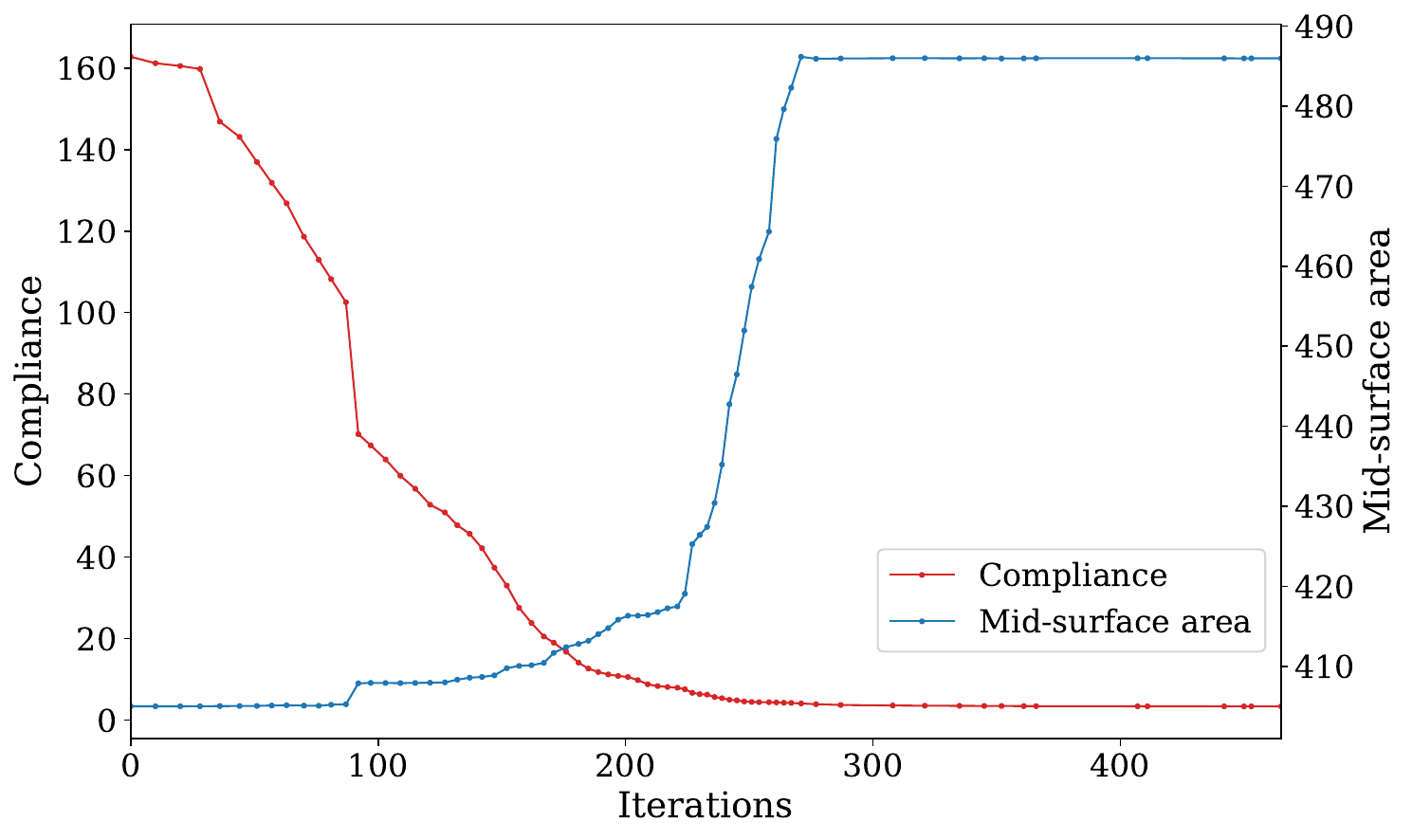}
		\caption{Compliance and mid-surface area variations during optimisation. The final compliance values are $2.541$ (left), $2.934$ (middle) and $3.356$ (right).}
	\end{subfigure}
	\\
	\begin{subfigure}{\linewidth}
		\includegraphics[width=0.25\linewidth]{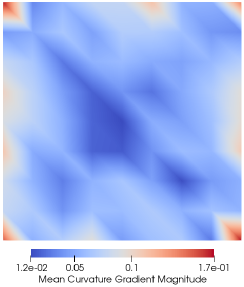}
		\hfill
		\includegraphics[width=0.25\linewidth]{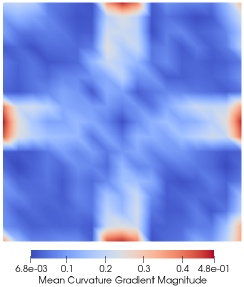}
		\hfill
		\includegraphics[width=0.25\linewidth]{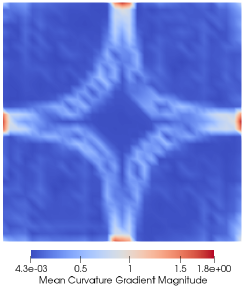}
		\caption{Mean curvature gradient magnitudes of the optimised shapes with $64$ (left), $256$ (middle) and $1024$ (right) shell elements.}
		\label{fig:shellOptCurvatureGrad}
	\end{subfigure}
	\caption{Shape optimisation of the square thin-shell roof under uniform loading with different mesh densities.}
	\label{fig:shellOptCompare}
\end{figure*}

A neural network with two hidden layers is used. Each hidden layer has five nodes, one layer with a periodic activation and the other with the ReLU activation. The number of network parameters is $51$. A coarse shell mesh contains $64$ elements, as shown in the left column of Figure~\ref{fig:shellOptCompare}. The structural compliance of the initial thin-shell is $130.444$, which is reduced to $2.541$ after optimisation.

To investigate the influence of the mid-surface discretisation, the mesh is subdivided into $256$ and $1024$ shell elements. The optimised thin-shells have similar compliances and shapes across different mesh densities, see Figure~\ref{fig:shellOptCompare}, and they were obtained with the same neural network architecture, i.e., the same number of network parameters as well. The results demonstrate that using NRep can substantially reduce the number of design variables in optimisation while preserving geometric details. In addition, the gradient magnitudes of mean curvature for the optimised shapes are plotted in Figure~\ref{fig:shellOptCurvatureGrad}, which show smooth and gradual variations in the surface curvature without noticeable abrupt changes, indicating the absence of significant oscillatory features.

\begin{figure*}
	\centering
		\includegraphics[width=0.32\linewidth]{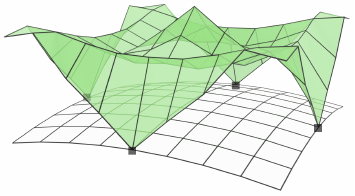}
		\hfill
		\includegraphics[width=0.32\linewidth]{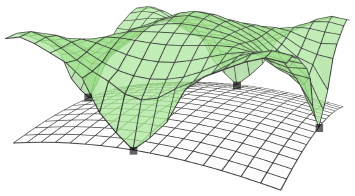}
		\hfill
		\includegraphics[width=0.32\linewidth]{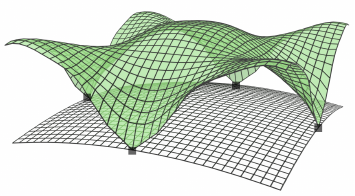}
	\caption{Optimised shapes under uniform loading obtained with the free-form deformation method for $64$ (left), $256$ (middle) and $1024$ (right) shell elements.}
	\label{fig:shellOptFFD}
\end{figure*}

To verify the optimisation results, the thin-shell roofs are also optimised using the free-form deformation (FFD) method~\cite{xiao2022free}, as shown in Figure~\ref{fig:shellOptFFD}. For the meshes consisting of $256$ and $1024$ shell elements, the FFD method yields the shapes that are consistent with those obtained using NRep. However, for the coarse mesh with $64$ shell elements, the FFD method converges to a local optimum characterised by a wavy surface, whereas the proposed method still captures the overall optimal shape.

\subsection{Thin-shell roof design under varied conditions}

\begin{figure*}
	\centering
	\begin{subfigure}{\linewidth}
		\includegraphics[width=0.3\linewidth]{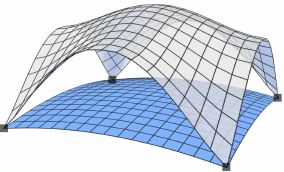}
		\hfill
		\includegraphics[width=0.3\linewidth]{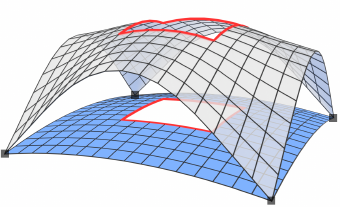}
		\hfill
		\includegraphics[width=0.3\linewidth]{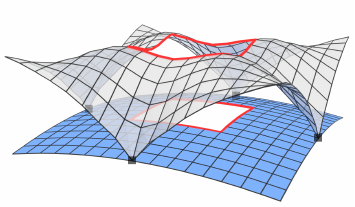}
		\caption{Uniform loading and fixed corners (left and middle) or mid-edge points (right)}
		\label{fig:uniformLoadExamples}
	\end{subfigure}
	\begin{subfigure}{\linewidth}
		\includegraphics[width=0.3\linewidth]{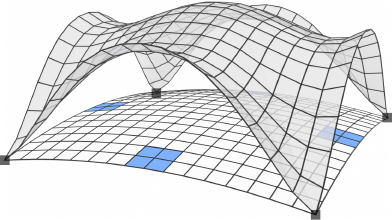}
		\hfill
		\includegraphics[width=0.3\linewidth]{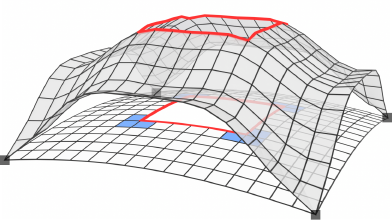}
		\hfill
		\includegraphics[width=0.3\linewidth]{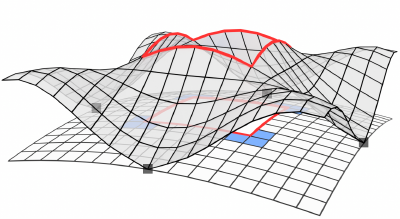}
		\caption{Regional loading with fixed corners (left and middle) or mid-edge points (right)}
		\label{fig:regionalLoadExamples}
	\end{subfigure}
	\caption{Shape optimisation with different boundary conditions and initial thin-shell geometries. The blue regions indicate the areas on which loading is applied and the red polygons highlight the central openings.}
	\label{fig:loadAndTopologyExamples}
\end{figure*}

The proposed approach is further applied to shape optimisation of thin-shells with varying topologies (with or without central openings) and boundary conditions (loads and supports), using a single fixed neural network architecture. Figure~\ref{fig:uniformLoadExamples} presents the optimisation results under uniform loading for three configurations: (i) a square thin-shell fixed at the four corners, (ii) a square thin-shell with a central square opening and corner supports, and (iii) the same thin-shell with supports relocated to the four mid-edge points. In case (i) the optimised shape exhibits a symmetric and dome-like form with smooth curvature, as expected for corner-supported thin-shells under uniform loading. Introducing a central opening (case (ii)) preserves the overall dome shape, while maintaining smooth curvature transitions and a regular boundary around the opening, indicating that the NRep captures the topology without generating spurious oscillations. In case (iii), relocating the supports to the mid-edge points substantially changes the load path in the thin-shell, resulting in a shape resembling the result without the opening (cf. Figure~\ref{fig:shellOptCompare}), with diagonal ridging patterns and more significant curvature variations.

For regional loading, the corresponding optimised shapes are shown in Figure~\ref{fig:regionalLoadExamples}. Compared with uniform loading, regional loading induces localised uplift in the loaded areas. In case (i), the surface exhibits a concentrated uplift in the mid-edge regions, with doubly curved boundaries around the loaded areas. For the thin-shell with a central opening (case (ii)), where loads are applied around the opening corners, the region surrounding the opening is uplifted relative to the shell edges, demonstrating the localised impact of the loads and the opening on the shape. In case (iii), the uplift around the opening remains the dominant feature but the overall shape changes due to the modified load path. In all these cases, the presence of the central opening does not introduce artefacts; instead, the NRep generates smooth boundary curves and continuous curvature transitions across the opening as well as in the highly curved regions where the loads are applied.

From all the cases in Figure~\ref{fig:loadAndTopologyExamples}, it is clear that even under complex combinations of supports, openings and non-uniform loads, the NRep generates stable and smooth geometries, effectively capturing both global structural behaviour and local geometric features without requiring modifications to the neural network architecture for each case.

\begin{figure*}
	\centering
		\includegraphics[width=0.3\linewidth]{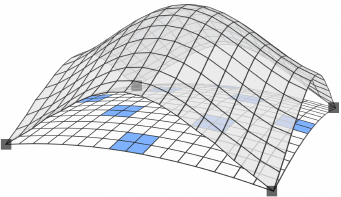}
		\hfill
		\includegraphics[width=0.3\linewidth]{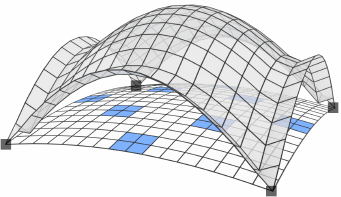}
		\hfill
		\includegraphics[width=0.3\linewidth]{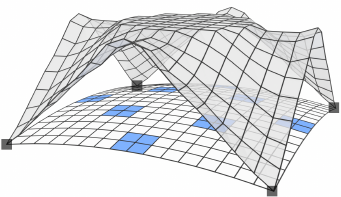}
	\caption{Optimised shapes obtained with NRep, considering each hidden layer containing $2$ nodes (left), $5$ nodes (middle) and $10$ nodes (right), with corresponding compliance $22.17$, $14.26$ and $9.54$, respectively. The blue regions indicate the areas on which loading is applied.}
	\label{fig:shellOptRegionalLoadNodes}
\end{figure*}

\begin{figure*}
	\centering
		\includegraphics[width=0.22\linewidth]{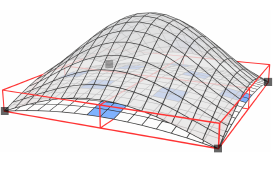}
		\hfill
		\includegraphics[width=0.22\linewidth]{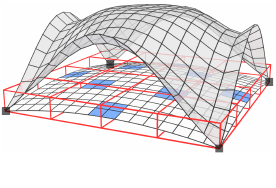}
		\hfill
		\includegraphics[width=0.22\linewidth]{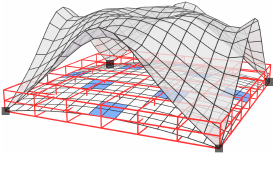}
		\hfill
		\includegraphics[width=0.3\linewidth]{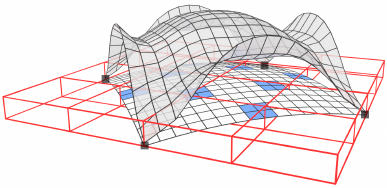}
	\caption{Optimised shapes obtained with the FFD method, considering the degrees of control volumes $2 \times 2 \times 1$, $4 \times 4 \times 1$ and $6 \times 6 \times 2$ and a rotated control volume. The blue regions indicate the areas on which loading is applied and the red polygons represent the initial FFD control volume.}
	\label{fig:shellOptRegionalLoadFFD}
\end{figure*}

\subsection{Parameter analysis and comparison}
To further assess the proposed method, the influence of neural network hyperparameters on the optimisation results of thin-shell roofs is investigated. The thin-shell roof is fixed at its corners and subjected to regional loading, as illustrated in Figure~\ref{fig:shellOptRegionalLoadNodes}. The initial geometry and material properties are identical to those used in the preceding examples. In the NRep framework, a two-layer MLP is adopted, with one hidden layer using a periodic activation (the period coefficient defaults to one if not explicitly stated) and the other using the ReLU activation.

First, the influence of the number of design variables is examined. In the proposed method, the number of design variables is determined by the number of nodes in the hidden layers. Three network configurations are considered, with each of the two hidden layers containing $2$, $5$ and $10$ nodes, resulting in $15$, $51$ and $151$ network parameters, respectively. The corresponding optimisation results are shown in Figure~\ref{fig:shellOptRegionalLoadNodes}. When each hidden layer contains only two nodes, the optimisation converges to an optimum that fails to capture the highly curved shape. As the number of nodes increases, the admissible design space is enlarged, allowing the optimisation to recover finer geometric features.

The optimisation results are further compared with those obtained using the FFD method. The degrees of the FFD control volumes are chosen as $2 \times 2 \times 1$, $4 \times 4 \times 1$ and $6 \times 6 \times 2$, corresponding to $18$, $50$ and $147$ design variables, respectively, such that the number of design variables is comparable with those in the NRep method. The resulting optimised shapes are presented in Figure~\ref{fig:shellOptRegionalLoadFFD}. As expected, higher-degree control volumes are capable of representing more geometric details, whereas control volumes with lower degrees fail to approach the optimal shape. Overall, the NRep and FFD methods exhibit similar trends with respect to the number of design variables. However, the optimisation results obtained with FFD are also influenced by the orientation of the thin-shell within the control volume. As shown in the last column of Figure~\ref{fig:shellOptRegionalLoadFFD}, when the thin-shell is rotated inside the control volume, the optimisation converges to a different optimal solution, since the boundary conditions are not straightforward to impose explicitly in the FFD method when the fixed positions do not coincide with the corners of the control volume. In contrast, the proposed NRep method is independent of the shell orientation and therefore avoids this limitation.

\begin{figure}
	\centering
	\begin{subfigure}{\linewidth}
		\includegraphics[width=0.55\linewidth]{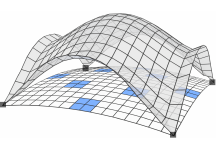}
		\hfill
		\includegraphics[width=0.4\linewidth]{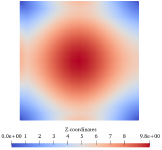}
		\caption{Period coefficient $\omega = 0.5$}
	\end{subfigure}
	\\
	\begin{subfigure}{\linewidth}
		\includegraphics[width=0.55\linewidth]{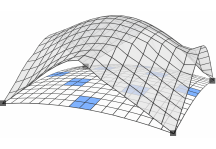}
		\hfill
		\includegraphics[width=0.4\linewidth]{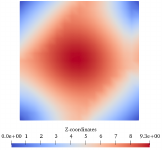}
		\caption{Period coefficient $\omega = 2$}
	\end{subfigure}
	\caption{Optimised shapes obtained using different periodic coefficients in activation functions and their $z$-coordinate contours. The blue regions in the initial mid-surface indicate the areas on which loading is applied.}
	\label{fig:shellOptRegionalLoadPeriod}
\end{figure}

In addition, the period coefficient in the periodic activation function (cf.~\eqref{eq:periodicActivation}) affects the surface approximation and may consequently influence the optimisation results. Two other values of the period coefficient are considered, $\omega = \{0.5, \, 2\}$. The optimisation results for $\omega = 0.5$ and $\omega = 2$ are shown in Figure~\ref{fig:shellOptRegionalLoadPeriod}. Since a larger period coefficient corresponds to approximation with higher frequency, the optimised shape obtained with $\omega = 2$ converges quickly to a local optimum and is less symmetric than those obtained using smaller values of $\omega$, as indicated by the contours of $z$-coordinates of the optimised shapes, whereas the solutions corresponding to $\omega = 0.5$ and $\omega = 1$ remain smooth and nearly symmetric. This indicates that moderate values of the period coefficient appear sufficient to capture the essential geometric features, while optimisation with excessively large values may be distracted by low-frequency geometric features.

\subsection{Application to lattice-skin geometry generation}
\begin{figure*}
	\centering
	\includegraphics[width=0.95\linewidth]{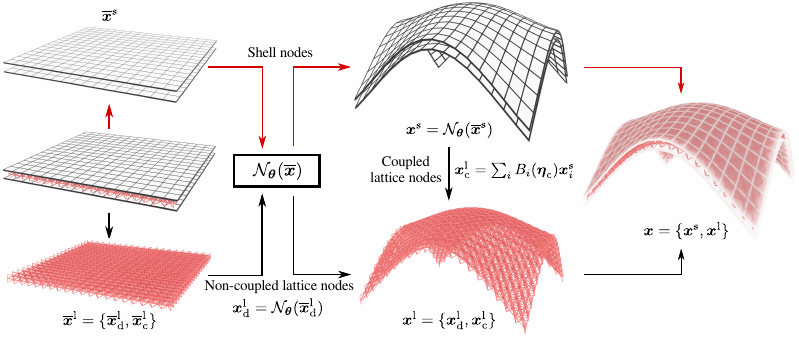}
	\caption{Lattice-skin geometry generated with NRep.}
	\label{fig:latticeSkinGeometry}
\end{figure*}

We next illustrate how the NRep of an optimised shell can be used in downstream structural design tasks. A lattice-skin structure is considered, comprising two thin-shells and an internal lattice, with the corners of the thin-shells fixed (Figure~\ref{fig:latticeSkinGeometry}). The projected area of the thin-shells is $20 \times 20$. The shell thickness is $0.1$ and the strut diameter is $0.04$. The lattice consists of two layers of body-centred unit cells with a total height of $1$. The Young's modulus and Poisson's ratio of the material are set to $7 \times 10^7$ and $0.35$, respectively.

The lower thin-shell of the lattice-skin is obtained from the optimised shape under uniform loading (see Figure~\ref{fig:uniformLoadExamples}). Figure~\ref{fig:latticeSkinGeometry} illustrates the procedure to generate the lattice-skin geometry with NRep. The upper thin-shell is generated by offsetting the lower thin-shell by $1$ in the vertical direction. Parametric coordinates of thin-shells and lattice nodes are computed by linearly scaling the physical coordinates into a unit cube $[0, 1]\times[0, 1]\times[0, 1]$. An NRep $\mathcal{N}_{\vec{\theta}}$ is then defined to map parametric coordinates of thin-shells to physical coordinates $\vec{x}^{\text{s}}\in\mathbb{R}^3$. For the lattice, positions of non-coupled lattice nodes, $\vec{x}_{\text{d}}^{\text{l}}$, are obtained directly using the same NRep $\mathcal{N}_{\vec{\theta}}$, whereas positions of coupled lattice nodes are determined from the corresponding coupling locations $\vec{\eta}_{\text{c}}$ on the thin-shells using shape functions $B_i$ and the coordinates of the associated shell nodes $\vec{x}_i^{\text{s}}$, i.e., $\vec{x}_{\text{c}}^{\text{l}} = \sum_i B_i(\vec{\eta}_\text{c})\vec{x}_i^{\text{s}}$, enforcing the displacement coupling between the thin-shells and the lattice~\cite{xiao2022infill}.

\begin{figure}
	\centering
	\begin{subfigure}{0.8\linewidth}
		\includegraphics[width=\linewidth]{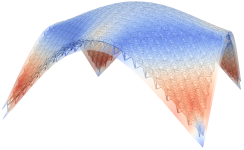}
		\caption{Isometric view.}
		\label{fig:latticeSkinDisp}
	\end{subfigure}
	\begin{subfigure}{0.8\linewidth}
		\includegraphics[width=\linewidth]{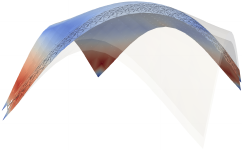}
		\caption{Clipped view.}
		\label{fig:latticeSkinDispCut}
	\end{subfigure}
	\caption{Structural analysis of the generated lattice-skin structure under mid-edge regional loads and isocontours of the displacement field.}
\end{figure}

The coupling information between the thin-shells and the lattice is preserved throughout the geometry generation. Consequently, the lattice-skin structure can be analysed immediately using the numerical technique proposed in~\cite{xiao2019interrogation} after geometry generation, therefore integrating the lattice-skin design and analysis. Figure~\ref{fig:latticeSkinDisp} shows the displacement field of the lattice-skin under a regional load of magnitude $10$ in the vicinity of the mid-edge points (cf. Figure~\ref{fig:regionalLoadExamples}). The clipped view in Figure~\ref{fig:latticeSkinDispCut} also confirms that the lattice generated with NRep conforms to the thin-shell geometry.

\section{Conclusion} \label{sec:conclusion}
We have introduced a novel shape optimisation framework in which the mid-surface of a thin-shell is parameterised using a neural parametric representation (NRep). The NRep defines a smooth mapping from parametric coordinates to three-dimensional physical space, therefore reformulating the shape optimisation problem in terms of neural network parameters. The NRep is flexible for geometric modelling and can be tuned through the choice of network architecture and activation functions; in particular, sinusoidal activations enable the faithful representation of fine geometric features with a relatively small number of network parameters. Consequently, thin-shells with fine finite element discretisations can be optimised using a low-dimensional design space, which reduces the dimensionality of the optimisation problem. The proposed approach integrates naturally with gradient-based constrained optimisation algorithms through automatic differentiation of the neural network. Numerical examples involving varying initial geometries, load cases and boundary conditions demonstrate that the proposed approach yields physically plausible optimised shapes and exhibits potential for generating more complex lattice-skin geometries.

While the examples presented here adopt fully connected multi-layer perceptrons, it is possible to use alternative neural network architectures, which may further improve approximation properties and optimisation convergence. Moreover, the NRep itself can be used to explore multiple distinct optima for the same problem, for instance, by heuristic deactivation or dropout of network parameters during optimisation to steer the search towards different regions of the design space. From a computational perspective, scalability to large-scale problems can be achieved by exploiting GPU acceleration for sensitivity analysis and incorporating neural surrogate models for structural analysis. Finally, in light of recent advances in neural implicit representations for topology optimisation, a promising avenue for future research would be the combination of the two neural representations to enable a unified framework for concurrent shape and topology optimisation.

\backmatter

\bibliography{references}

\end{document}